\input amstex
\documentstyle {amsppt}
\pagewidth{12.5 cm}\pageheight{19 cm}\magnification\magstep1
\topmatter
\title Singular supports for character sheaves on a group compactification\endtitle
\author Xuhua He and George Lusztig\endauthor
\thanks{X. H. is is supported by NSF grant DMS-0111298. G.L. is supported in part by NSF grant DMS-0243345.}\endthanks
\address School of Mathematics, Institute for Advanced Study, Princeton, NJ 08540\endaddress
\email hugo\@math.ias.edu \endemail
\address Department of Mathematics, M.I.T., Cambridge, MA 02139\endaddress
\email gyuri\@math.mit.edu \endemail\subjclassyear{2000}

\abstract{Let $G$ be a semisimple adjoint group over $\bold C$ and
$\bar{G}$ be the De Concini-Procesi completion of $G$. In this
paper, we define a Lagrangian subvariety $\Lambda$ of the
cotangent bundle of $\bar{G}$ such that the singular support of
any character sheaf on $\bar{G}$ is contained in
$\Lambda$.}\endabstract
\endtopmatter
\document

\define\po{\text{\rm pos}}

\define\Lie{\text{\rm Lie}}
\define\nl{\newline}
\define\Ad{\text{\rm Ad}}
\redefine\i{^{-1}}

\define\cd{\Cal D}

\define\cn{\Cal N}
\define\co{\Cal O}
\define\cp{\Cal P}

\define\cz{\Cal Z}

\define\a{\alpha}

\define\G{\Gamma}
\redefine\d{\delta} \redefine\D{\Delta}

\define\p{\pi}

\define\x{\xi}

\define\T{\times}
\define\fg{\frak g}
\define\sub{\subset}
\define\sqc{\sqcup}
\define\m{\mapsto}
\define\bbq{\bar{\bold Q}_l}
\define\bsl{\backslash}
\redefine\L{\Lambda}
\define\bst{\bigstar}
\define\iy{\infty}
\define\WW{\bold W}
\define\II{\bold I}
\define\NN{\bold N}
\define\tB{\tilde B}
\define\tb{\tilde b}
\define\tX{\tilde X}
\define\em{\emptyset}
\define\bG{\bar G}
\define\si{\sim}
\define\dd{\bold d}
\define\un{\underline}
\redefine\ss{\bold s} \subhead 1.1\endsubhead In this paper all
algebraic varieties are assumed to be over a fixed algebraically
closed field of characteristic $0$.

If $X$ is a smooth variety, let $T^* X$ be the cotangent bundle of
$X$. For any morphism $\a: X \rightarrow Y$ of smooth varieties
and $x \in X$, we write $\a^*: T^*_{\a(x)}Y@>>>T^*_xX$ for the map
induced by $\a$. If moreover, $\a:X@>>>Y$ is a locally trivial
fibration with smooth connected fibres and $\L$ is a closed
Lagrangian subvariety of $T^*Y$, then let $\a^\bst(\L)=\cup_{x \in
X} \a^* (\L\cap T^*_{\a(x)}Y) \subset T^*X$. Then

(a) {\it $\a^\bst(\L)$ is a closed Lagrangian subvariety of
$T^*X$. Moreover, the set of irreducible components of $\L$ is
naturally in bijection with the set of irreducible components of
$\a^\bst(\L)$.}

Let $X,Y$ be smooth irreducible varieties and let $\a:X@>>>Y$ be a
principal $P$-bundle for a free action of a connected linear
algebraic group $P$ on $X$.

(b) {\it If $\L'$ is a closed Lagrangian subvariety of $T^*X$
stable under the $P$-action then $\L'=\a^\bst(\L)$ for a unique
closed Lagrangian subvariety $\L$ of $T^*Y$.}

Let $X$ be a smooth irreducible variety and let $i:Y@>>>X$ be the
inclusion of a locally closed smooth irreducible subvariety. Let
$\L$ be a closed Lagrangian subvariety of $T^*Y$. Let $i_\bst(\L)$
be the subset of $T^*X$ consisting of all $\x\in T^*_xX$ such that
$x\in Y$ and the image of $\x$ under the obvious surjective map
$T^*_xX@>>>T^*_xY$ belongs to $\L\cap T^*_xY$. Note that

(c) {\it $i_\bst(\L)$ is a locally closed Lagrangian subvariety of
$T^*X$. Moreover, the set of irreducible components of
$i_\bst(\L)$ is naturally in bijection with the set of irreducible
components of $\L$.}

For an algebraic variety $X$ we write $\cd(X)$ for the bounded
derived category of constructible $\bbq$-sheaves on $X$ where $l$
is a fixed prime number. For $X$ smooth and $C\in\cd(X)$, we
denote by $SS(X)$ the singular support of $C$ (a closed Lagrangian
subvariety of $T^*X$). Let $A$ be a connected linear algebraic
group acting on a smooth variety $X$ and let $B$ be a connected
subgroup of $A$. Let $\mu_A: T^*X@>>>\Lie(A)^*$ be the moment map
of the $A$-action on $X$. Consider the diagram $X@<pr_1<<A\T
X@>pr_2>>A\T_BX@>\p>>X$ where $B$ acts on $A\T X$ by
$b:(a,x)\m(ab\i,bx)$, $A\T_B X$ is the quotient space and
$\p(a,x)=ax$. Then for any $B$-equivariant perverse sheaf $C$ on
$X$ there is a well defined perverse sheaf $C'$ on $A\T_BX$ such
that $pr_2^*C'=pr_1^*C$ up to a shift. We set
$\G^A_B(C)=\p_*C'\in\cd(X)$. By \cite{MV, 1.2} we have

(d) $SS(\G^A_B(C))\sub\overline{A\cdot SS(C)}$. \nl On the other
hand, we have

(e) $SS(C)\sub\mu_B\i(0)$. \nl Indeed, if $p_1:B\T X@>>>X$ is the
action and $p_2:B\T X@>>>X$ is the second projection we have
$p_1^*(C)=p_2^*(C)$. Hence $SS(p_1^*(C))=SS(p_2^*(C))$. Using
\cite{KS, 4.1.2} we can rewrite this as
$p_1^\bst(SS(C))=p_2^\bst(SS(C))$. Hence if $x\in X$ and $\x\in
T^*_xX\cap SS(C)$ then the image of $\x$ under the map
$T^*_xX@>>>T^*_1(B)$ induced by $B@>>>X,b\m bx$ is $0$. This
proves (e).

\subhead 1.2\endsubhead Let $G$ be a connected reductive algebraic
group. Let $\fg=\Lie(G)$. Let $\cn$ be the variety of nilpotent
elements in $\fg^*$. Let $B$ be a Borel subgroup of $G$. Let $K$
be a closed connected subgroup of $G$. Then $B_K=B\cap K$ is a
parabolic subgroup of $K$. Assume that $G$ acts on a smooth
variety $X$. Let $C$ be a $B_K$-equivariant perverse sheaf on $X$;
assume also that there exists a finite covering $a:\tB@>>>B$ such
that $C$ is $\tB$-equivariant for the $\tB$-action $\tb:x\m
a(\tb)x$ on $X$. By 1.1(e) we have $\mu_{\tB}(SS(C))=0$. Since
$\Lie(\tB)=\Lie(B)$ we then have $\mu_B(SS(C))=0$. It follows that
$\mu_G(SS(C))$ is contained in the kernel of the obvious map
$\fg^*@>>>\Lie(B)^*$ hence is contained in $\cn$. Since $\cn$ is
stable under the coadjoint action we have $\mu_G(K\cdot
SS(C))=K\mu_G(SS(C))\sub\cn$. Using this together with 1.1(d) and
the fact that $\mu_G\i(\cn)$ is closed in $T^*X$ we see that
$SS(\G^K_{B_K}(C))\sub\mu_G\i(\cn)$. Applying 1.1(e) to
$\G^K_{B_K}(C),K$ instead of $C,B$ we see that
$SS(\G^K_{B_K}(C))\sub\mu_K\i(0)=\mu_G\i(\Lie(K)^\perp)$ where
$\Lie(K)^\perp\sub\fg^*$ is the annihilator of $\Lie(K)\sub\fg$.
Thus we have
$$SS(\G^K_{B_K}(C))\sub\mu_G\i(\Lie(K)^\perp\cap\cn).\tag a$$

\subhead 1.3\endsubhead We now replace $G,\fg,B,K,B_K,X,C$ by $G\T
G,\fg\T\fg,B\T B,G_\D,B_\D,X',C'$ where $G_\D=\{(g,g')\in G\T
G;g=g'\},B_\D=\{(g,g')\in B\T B;g=g'\}$, $X'$ is a smooth variety
with a given action of $G\T G$ and $C'$ is a $B_\D$-equivariant
perverse sheaf on $X'$; we assume that there exists a finite
covering $a':\tB'@>>>B\T B$ such that $C'$ is $\tB'$-equivariant
for the $\tB'$-action $\tb':x'\m a'(\tb')x'$ on $X'$. We have the
following special case of 1.2(a):
$$SS(\G^{G_\D}_{B_\D}(C'))\sub\mu_{G\T G}\i(\cn^-)\tag a$$
where
$$\cn^-=\{(f,f')\in\fg^*\T\fg^*;f+f'=0,f,f'\text{ nilpotent }\}.$$

\subhead 1.4\endsubhead Let $\WW$ be the Weyl group of $G$ and let
$\II$ be the set of simple reflections in $\WW$. Let $G'$ be a
possibly disconnected algebraic group with identity component $G$
and with a given connected component $D$. Now $G\T G$ acts
transitively on $D$ by $(g_1,g_2):g\m g_1gg_2\i$. Hence the moment
map $\mu_{G\T G}:T^*D@>>>\fg^*\T\fg^*$ is well defined. In
\cite{L1, 4.5} a class of perverse sheaves (called character
sheaves) on $D$ is introduced. These appear as constituents of
some perverse cohomology sheaf of $\G^{G_\D}_{B_\D}(C')$ for some
$C'$ as in 1.3 (with $X'=D$). Hence from 1.3(a) we deduce:

(a) {\it If $K$ is a parabolic character sheaf on $D$ then
$SS(K)\sub\mu_{G\T G}\i(\cn^-)$.} \nl In the case where $G'=G=D$ a
statement close to (a) appears in \cite{MV, 2.8} (where it is
attributed to the second author) and in \cite{Gi}.

\subhead 1.5\endsubhead We preserve the setup of 1.4. For any
$J\sub\II$ let $\cp_J$ be the set of parabolic subgroups of $G$ of
type $J$. In particular $\cp_\em$ is the set of Borel subgroups of
$G$. For $J\sub\II$ let $\WW_J$ be the subgroup of $\WW$ generated
by $J$; let $\WW^J$ (resp. ${}^J\WW$) be the set of all $w\in\WW$
such that $w$ has minimal length among the elements in $\WW_Jw$
(resp. $w\WW_J$). Let $\d:\WW@>\si>>\WW$ be the isomorphism such
that $\d(\II)=\II$ and such that $J\sub\II,P\in\cp_J,g\in
D\implies gPg\i\in\cp_{\d(J)}$. Following \cite{L2, 8.18}, for
$J,J'\sub\II$ and $y\in{}^{J'}\WW\cap\WW^J$ such that
$\Ad(y)(\d(J))=J'$ we set
$$Z_{J,y,\d}=\{(P,P',gU_P);P\in\cp_J,P'\in\cp_{J'},g\in D,\po(P',gPg\i)=y\}.$$
Now $G\T G$ acts (transitively) on $Z_{J,y,\d}$ by
$$(g_1,g_2):(P,P',gU_P)\m(g_2Pg_2\i,g_1P'g_1\i,g_1gg_2\i).$$
Hence the moment map $\mu_{G\T G}:T^*Z_{J,y,\d}@>>>\fg^*\T\fg^*$
is well defined. In \cite{L2, Section 11} a class of perverse
sheaves (called parabolic character sheaves) on $Z_{J,y,\d}$ is
introduced. These appear as constituents of some perverse
cohomology sheaf of $\G^{G_\D}_{B_\D}(C')$ for some $C'$ as in 1.3
(with $X'=Z_{J,y,\d}$). Hence from 1.3(a) we deduce:

(a) {\it If $K$ is a parabolic character sheaf on $Z_{J,y,\d}$
then $SS(K)\sub\mu_{G\T G}\i(\cn^-)$.} \nl When $J=\II$, this
reduces to 1.4(a).

\subhead 1.6\endsubhead Assume that $G$ is adjoint. Let $\bG$ be
the De Concini-Procesi compactification of $G$. Then $G\T G$ acts
naturally on $\bG$ extending continuously the action
$(g_1,g_2):g\m g_1gg_2\i$ of $G\T G$ on $G$. Hence the moment map
$\mu_{G\T G}:T^*\bG@>>>\fg^*\T\fg^*$ is well defined. In \cite{L2}
a class of perverse sheaves (called parabolic character sheaves)
on $\bG$ is introduced. It has been shown by He \cite{H2} and by
Springer (unpublished) that any parabolic character sheaf on $\bG$
appears as a constituent of some perverse cohomology sheaf of
$\G^{G_\D}_{B_\D}(C')$ for some $C'$ as in 1.3 (with $X'=\bG$).
Hence from 1.3(a) we deduce:

(a) {\it If $K$ is a parabolic character sheaf on $\bG$ then
$SS(K)\sub\mu_{G\T G}\i(\cn^-)$.}

\subhead 1.7\endsubhead In the setup of 1.4 let $\L(D)=\mu_{G\T
G}\i(\cn^-)$. We want to describe the variety $\L(D)$. For $g\in
D$ let $I_g$ be the isotropy group at $g$ of the $G\T G$-action on
$D$ that is, $I_g=\{(g_1,g_2)\in G\T G;g_2=g\i g_1g\}$. We have
$\Lie(I_g)=\{(y_1,y_2)\in\fg\T\fg;y_2=\Ad(g)\i(y_1)\}$ and the
annihilator of $\Lie(I_g)$ in $\fg^*\T\fg^*$ is
$\Lie(I_g)^\perp=\{(z_1,z_2)\in\fg^*\T\fg^*;z_1+\Ad(g)(z_2)=0\}$.
This may be identified with the fibre of $T^*D$ at $g$. Then
$$\align&\L(D)=\{(g,z_1,z_2)\in D\T\fg^*\T\fg^*;z_1+\Ad(g)(z_2)=0,z_1+z_2=0,z_2\in\cn\}
\\&=\{(g,z,-z);g\in D,z\in\cn,\Ad(g)(z)=z\}=\sqc_\co X_\co\endalign$$
where $\co$ runs over the (finite) set of $\Ad(G)$-orbits on $\cn$
which are normalized by some element of $D$ and
$X_\co=\{(g,z,-z);g\in D,z\in\co,\Ad(g)(z)=z\}$. We pick
$\x\in\co$ and let $\cz'=\{h\in G';\Ad(h)\x=\x\}$, $\cz=\{h\in
G;\Ad(h)\x=\x\}$. Let $\un\cz'$ (resp. $\un\cz$) be the group of
connected components of $\cz'$ (resp. $\cz$). Then $\un\cz'$ is a
finite group and $\un\cz$ is a subgroup of $\un\cz'$. Let
$\un\cz_1$ be the set of connected components of $\cz'$ that are
contained in $D$. Then $\un\cz_1$ is a subset of $\cz'$; also,
$\un\cz$ acts on $\un\cz_1$ by conjugation inside $\un\cz'$. Let
$F_\co^D$ be the set of orbits of this action. Note that $F_\co^D$
is independent (up to unique isomorphism) of the choice of $\x$.

Let $\tX=\{(g,r)\in D\T G;r\i gr\in\cz'\}$. Then $\cz$ acts freely
on $\tX$ by $h:(g,r)\m(g,rh\i)$ and we have an isomorphism
$\cz\bsl\tX@>\si>>X_\co$, $(g,r)\m(g,\Ad(r)\x)$. By the change of
variable $(g,r)\m(g',r)$, $g'=r\i gr$, $\tX$ becomes
$\{(g',r);g'\in\cz'\cap D,r\in G\}$. In the new coordinates, the
free action of $\cz$ on $\tX$ is $h:(g',r)\m(hg'h\i,rh\i)$. We see
that $\tX$ is smooth of pure dimension $\dim(\cz\T G)$ and its
connected components are indexed naturally by $\un\cz_1$ (the
connected component containing $(g',r)$ is indexed by the image of
$g'$ in $\un\cz_1$). The action of $\cz$ on $\tX$ permutes the
connected components of $\tX$ according to the action of $\un\cz$
on $\un\cz_1$ considered above. We see that $X_\co=\cz\bsl\tX$ is
smooth of pure dimension $\dim G$ and its connected components are
indexed naturally by the set $F_\co^D$.

We see that $\L(D)$ can be partitioned into finitely many locally
closed, irreducible, smooth subvarieties of dimension $\dim G$,
indexed by the finite set $F(D):=\sqc_\co F^D_\co$. In particular,
$\L(D)$ has pure dimension $\dim G$. More precisely, one checks
that

(a) {\it $\L(D)$ is a closed Lagrangian subvariety of $T^*D$.}

\subhead 1.8\endsubhead In the setup of 1.5 we set
$\L(Z_{J,y,\d})=\mu_{G\T G}\i(\cn^-)$. We want to describe the
variety $\L(Z_{J,y,\d})$.

Following \cite{L2, 8.18} we consider the partition
$Z_{J,y,\d}=\sqc_\ss Z_{J,y,\d}^\ss$ where $Z_{J,y,\d}^\ss$ are
certain locally closed smooth irreducible $G_\D$-stable
subvarieties of $Z_{J,Y,\d}$ indexed by the elements $\ss$ of a
finite set $S(J,\Ad(y)\d)$ which is in canonical bijection with
${}^{J'}\WW$, see \cite{L1, 2.5}. Note that each $\ss$ is a
sequence $(J_n,J'_n,u_n)_{n\ge0}$ where $J_n,J'_n$ are subsets of
$\II$ such that $J_n,J'_n$ are independent of $n$ for large $n$
and $u_n\in\WW$ is $1$ for large $n$.

We wish to define a Lagrangian subvariety $\L(Z_{J,y,\d}^\ss)$ of
$T^*(Z_{J,y,\d}^\ss)$.

Assume first that $\ss$ is such that $J_n=J,J'_n=J',u_n=1$ for all
$n$. In this case we have $J=J'$. Let $P\in\cp_J$ and let $L$ be a
Levi subgroup of $P$. Then $\dd_\ss=\{g\in
D;gLg\i=L,\po(P,gPg\i)=y\}$ is a connected component of the
algebraic group $N_{G'}(L)$ with identity component $L$. Hence
$\L(\dd_\ss)\sub T^*\dd_\ss$ is defined as in 1.7. We have a
diagram $Z_{J,y,\d}^\ss@<\a<<G\T\dd_\ss@>pr_2>>\dd_\ss$ where
$\a(h,g)=(hPh\i,hPh\i,U_{hPh\i}hgh\i U_{hPh\i})$. Note that $\a$
is a principal $P$-bundle where $P$ acts on $G\T\dd_\ss$ by
$p:(h,g)=(hp\i,\bar pg\bar p\i)$ (we denote the canonical
homomorphism $P@>>>L$ by $\bar p$). Let
$\L'=pr_2^\star\L(\dd_\ss)\sub T^*(G\T\dd_\ss)$. By 1.1(a) and
1.7(a), $\L'$ is a closed Lagrangian subvariety of $T^*
(G\T\dd_\ss)$. It is clearly stable under the natural action of
$P$ on $T^*(G\T\dd_\ss)$ (since $\L(\dd_\ss)$ is $L$-stable). By
1.1(b) there is a unique Lagrangian subvariety $\L''$ of
$T^*(Z_{J,y,\d}^\ss)$ such that $\a^\bst\L''=\L'$. We set
$\L(Z_{J,y,\d}^\ss)=\L''$.

We now consider a general $\ss=(J_n,J'_n,u_n)_{n\in\NN}$. For any
$r\in\NN$ let $\ss_r=(J_n,J'_n,u_n)_{n\ge r}$, $y_r=u_{r-1}\i\dots
u_1\i u_0\i y$. Then $Z^{\ss_r}_{J_r,y_r,\d}$ is defined and we
have a canonical map
$f_r:Z^\ss_{J,y,\d}@>>>Z^{\ss_r}_{J_r,y_r,\d}$ (a composition of
affine space bundles, see \cite{L2, 8.20(a)}). Moreover for
sufficiently large $r$, $\ss_r,J_r,y_r,f_r$ are independent of
$r$; we write $\ss_\iy,J_\iy,y_\iy,f_\iy$ instead of
$\ss_r,J_r,y_r,f_r$. Note also that $\ss_\iy,J_\iy,y_\iy$ are of
the type considered earlier, so that
$\L(Z^{\ss_\iy}_{J_\iy,y_\iy,\d})$ is defined as above. We set
$\L(Z_{J,y,\d}^\ss)=f_\iy^\bst(\L(Z^{\ss_\iy}_{J_\iy,y_\iy,\d}))$.

We now define
$$\L'(Z_{J,y,\d})=\sqc_{\ss\in S(J,\Ad(y)\d)}(i_\ss)_\bst(\L(Z_{J,y,\d}^\ss))$$
where $i_\ss:Z_{J,y,\d}^\ss@>>>Z_{J,y,\d}$ is the inclusion. From
1.1(c) we see that $\L'(Z_{J,y,\d})$ is a finite union of locally
closed Lagrangian subvarieties of $T^*(Z_{J,y,\d})$.

We state the following result: \proclaim{Proposition 1.9} We have
$\L(Z_{J,y,\d})=\L'(Z_{J,y,\d})$. In particular, $\L'(Z_{J,y,\d})$
is closed in $T^*(Z_{J,y,\d})$ and $\L(Z_{J,y,\d})$ is a
Lagrangian subvariety of $T^*(Z_{J,y,\d})$.
\endproclaim

Let $x \in Z_{J, y, \d}^\ss$. We will show that $\L(Z_{J, y, \d})
\cap T^*_x(Z_{J, y, \d})=\L'(Z_{J, y, \d}) \cap T^*_x(Z_{J, y,
\d})$.

We identify $\fg$ with $\fg^*$ via a $G$-invariant symmetric
bilinear form. Choose an element $g_0$ in $D$ that normalizes $B$
and a maximal torus $T$ of $B$. Let $P_J$ be the unique element in
$\cp_J$ that contains $B$. For $w \in W$, we choose a
representative $\dot w$ of $w$ in $N(T)$. Set $h_{J, y, \d}=(P_J,
{}^{\dot y \i} P_{J'}, U_{^{\dot y \i} P_{J'}} g_0 U_{P_J}) \in
Z_{J, y, \d}$. The element $\ss \in S(J,\Ad(y)\d)$ corresponds to
an element $w' \in {}^{J'} W$ under the bijection in \cite{L1,
2.5}. Set $w=(w')\i y \in w \in W^{\d(J)}$. By \cite{H1, 1.10},
$$Z_{J, y, \d}^\ss=G_{\D} (P_{J_{\infty}} \dot w, 1) \cdot h_{J, y,
\d}=G_{\D} (L_{J_{\infty}} \dot w, 1) \cdot h_{J, y, \d}.$$

Note that $\L(Z_{J,y,\d})$ and $\L'(Z_{J,y,\d})$ are stable under
the action of $G_{\D}$. Then we may assume that $x=(l_1 \dot w, 1)
\cdot h_{J, y, \d}$ for some $l_1 \in L_{J_{\infty}}$. By
\cite{H1, 1.10(2)}, $(P_{J_{\infty}})_{\D} \cdot (L_{J_{\infty}},
1) x=(P_{J_{\infty}}, 1) \cdot x$, where
$(P_{J_{\infty}})_{\D}=\{(g, g') \in P_{J_{\infty}} \times
P_{J_{\infty}}; g=g'\}$. For $u \in U_{P_{J_{\infty}}}$, $(u, 1)
\cdot x=(u' l, u') \cdot x$ for some $u' \in U_{P_{J_{\infty}}}$
and $l \in L_{J_{\infty}}$. Note that $(U_{P_{J_{\infty}}},
U_{P_{J_{\infty}}}) \cdot (l, 1) \cdot x=(U_{P_{J_{\infty}}} l, 1)
\cdot x$ and $(P_{J_{\infty}}, 1) \cdot x \cong U_{P_{J_{\infty}}}
\times (L_{J_{\infty}}, 1) \cdot x$. Thus $(l, 1) \cdot x=x$ and
$(U_{P_{J_{\infty}}}, 1) \cdot x \subset (P_{J_{\infty}})_{\D}
\cdot x$. Therefore, for $(f, -f) \in \mu_{G \times G}(T^*_x(Z_{J,
y, \d}))$, we have that $(f, -f)(\Lie(U_{P_{J_{\infty}}}), 0)=(f,
-f)(\Lie(P_{J_{\infty}})_{\D})=0$, i. e., $f \in
\Lie(P_{J_{\infty}})$. In particular, $f$ is nilpotent if and only
if the image of $f$ under $\Lie(P_{J_{\infty}}) \rightarrow
\Lie(P_{J_{\infty}})/\Lie(U_{P_{J_{\infty}}}) \cong
\Lie(L_{J_{\infty}})$ is nilpotent.

Hence $\mu_{G \times G}(\L(Z_{J, y, \d}) \cap T^*_x(Z_{J, y,
\d}))$ consists elements of the form $(u+l, -u-l)$ with $u \in
\Lie(U_{P_{J_{\infty}}})$, $l$ nilpotent in $\Lie(L_{J_{\infty}})$
and $(u, -u)I_x=(l, -l)I_x=0$, where $I_x$ is the Lie subalgebra
of the isotropic subgroup of $G \times G$ at point $x$.

Denote by $N_x$ the stalk at point $x$ of the conormal bundle
$N^*_{Z_{J, y, \d}^\ss}(Z_{J, y, \d})$. Since $Z_{J, y,
\d}^\ss=G_{\D}(P_{J_{\infty}}, 1) h_{J, y, \d}$, we have $$\mu_{G
\times G}(N_x)=\{(u, -u) ; u \in \Lie(U_{P_{J_{\infty}}}), (u,
-u)I_x=0\}.$$

Let $p_x: T^*_x(Z_{J, y, \d}) \rightarrow T^*_x(Z_{J, y, \d}^\ss)
\cong T^*_x(Z_{J, y, \d})/N_x$ be the obvious surjective map. Then
$\L(Z_{J, y, \d}) \cap T^*_x(Z_{J, y, \d})=p_x \i
\bigl(p_x(\L(Z_{J, y, \d}) \cap T^*_x(Z_{J, y, \d})) \bigr)$. Note
that $$I_x=\{(u_1+\Ad(l_1 \dot w g_0) l, u_2+l); u_1 \in \Ad(l_1
\dot w \dot y \i) \Lie(U_{P_{J'}}), u_2 \in \Lie(U_{P_J}), l \in
\Lie(L_J)\}.$$ Thus for $l \in \Lie(L_{J_{\infty}})$, $(l, -l)
I_x=0$ if and only if $\Ad(l_1 \dot w g_0) l=l$.

We identify $T^*_x(Z_{J, y, \d})$ with $\mu_{G \times G}
\bigl(T^*_x(Z_{J, y, \d})\bigr) \subset \fg \times \fg$ and regard
$T^*_x(Z_{J, y, \d}^\ss)$ as a subspace of $(\fg \times \fg)/N_x$.
Set $M=\{(l, -l); l \in \Lie(L_{J_{\infty}}), \Ad(l_1 \dot w g_0)
l=l\} \subset \fg \times \fg$. Then

(1) $p_x(\L(Z_{J, y, \d}) \cap T^*_x(Z_{J, y, \d}))=(M+N_x)/N_x$.

Now consider the commuting diagram
$$\CD
(P_{J_{\infty}})_{\D} (L_{J_{\infty}} \dot w, 1) \cdot h_{J, y,
\d} @>{i_1}>>Z_{J, y, \d}^\ss \cr @V{f'_{\infty}}VV
@V{f_{\infty}}VV \cr (P_{J_{\infty}})_{\D} (L_{J_{\infty}} \dot w,
1) \cdot h_{J_{\infty}, y_{\infty}, \d} @>{i_2}>>Z_{J_{\infty},
y_{\infty}, \d}^{\ss_{\infty}}, \endCD$$ where $i_1, i_2$ are
inclusions and $f'_{\infty}$ is the restriction of $f_{\infty}$.
Let $x'=f_{\infty}(x) \in Z_{J_{\infty}, y_{\infty},
\d}^{\ss_{\infty}}$. Since $f'_{\infty}$ is $P_{J_{\infty}} \times
P_{J_{\infty}}$-invariant, we have the following commuting diagram

$$\CD
T^*_{x'}\bigl((P_{J_{\infty}})_{\D} (L_{J_{\infty}} \dot w, 1)
\cdot h_{J_{\infty}, y_{\infty}, \d} \bigr)
@>{{}^1\mu_{P_{J_{\infty}} \times
P_{J_{\infty}}}}>>\Lie(P_{J_{\infty}})^* \times
\Lie(P_{J_{\infty}})^* \cr @V{(f'_{\infty})^*}VV @V{id}VV \cr
T^*_x \bigl((P_{J_{\infty}})_{\D} (L_{J_{\infty}} \dot w, 1) \cdot
h_{J, y, \d} \bigr) @>{{}^2\mu_{P_{J_{\infty}} \times
P_{J_{\infty}}}}>>\Lie(P_{J_{\infty}})^* \times
\Lie(P_{J_{\infty}})^*,
\endCD$$ where ${}^1\mu_{P_{J_{\infty}} \times
P_{J_{\infty}}}$ and ${}^2\mu_{P_{J_{\infty}} \times
P_{J_{\infty}}}$ are the moment maps. Since the actions of
$P_{J_{\infty}} \times P_{J_{\infty}}$ on $(P_{J_{\infty}})_{\D}
(L_{J_{\infty}} \dot w, 1) \cdot h_{J, y, \d}$ and
$(P_{J_{\infty}})_{\D} (L_{J_{\infty}} \dot w, 1) \cdot
h_{J_{\infty}, y_{\infty}, \d}$ are transitive,
${}^1\mu_{P_{J_{\infty}} \times P_{J_{\infty}}}$ and
${}^2\mu_{P_{J_{\infty}} \times P_{J_{\infty}}}$ are injective.

Set $\L_{x'}(Z_{J_{\infty}, y_{\infty},
\d}^{\ss_{\infty}})=\L(Z_{J_{\infty}, y_{\infty},
\d}^{\ss_{\infty}}) \cap T^*_{x'}(Z_{J_{\infty}, y_{\infty},
\d}^{\ss_{\infty}})$. Then $${}^2\mu_{P_{J_{\infty}} \times
P_{J_{\infty}}} \bigl(i^*_2(\L_{x'}(Z_{J_{\infty}, y_{\infty},
\d}^{\ss_{\infty}}) \bigr)={}^1\mu_{P_{J_{\infty}} \times
P_{J_{\infty}}} \bigl((f'_{\infty})^* i^*_2(\L_{x'}(Z_{J_{\infty},
y_{\infty}, \d}^{\ss_{\infty}}) \bigr)=M.$$ Here we identify
$\Lie(P_{J_{\infty}})^*$ with $\Lie(P^-_{J_{\infty}})$ via the
symmetric bilinear form. Moreover, $i_1^*
f^*_{\infty}=(f'_{\infty})^* i_2^*$ maps $\L_{x'}(Z_{J_{\infty},
y_{\infty}, \d}^{\ss_{\infty}})$ bijectively onto its image.
Therefore ${}^1\mu_{P_{J_{\infty}} \times P_{J_{\infty}}}
\bigl(i^*_1(\L(Z_{J, y, \d}^{\ss}) \cap T^*_x(Z_{J, y, \d}^{\ss}))
\bigr)=M$ and $i_1^*$ maps $\L(Z_{J, y, \d}^{\ss}) \cap
T^*_x(Z_{J, y, \d}^{\ss})$ bijectively onto its image. In other
words,

(2) $\L(Z_{J, y, \d}^{\ss}) \cap T^*_x(Z_{J, y,
\d}^{\ss})=(M+N_x)/N_x$.

Combining (1) and (2), $p_x(\L(Z_{J, y, \d}) \cap T^*_x(Z_{J, y,
\d}))=\L(Z_{J, y, \d}^\ss) \cap T^*_x(Z_{J, y, \d}^\ss)$. The
proposition is proved.

\proclaim{Corollary 1.10} The set of irreducible components of
$\L(Z_{J,y,\d})$ is in natural bijection with $\sqc_{\ss\in
S(J,\Ad(y)\d)}F(\dd_{\ss_\iy})$ (notation of 1.7, 1.8).
\endproclaim

\subhead 1.11\endsubhead In the setup of 1.6 let $\L(\bG)=\mu_{G\T
G}\i(\cn^-)$. We want to describe the variety $\L(Z_{J,y,\d})$. As
in \cite{L2, 12.3}, we have $\bG=\sqc_{J\sub\II}G_J$ where $G_J$
are the various $G\T G$-orbits in $\bG$; moreover we may identify
$G_J=T_J\bsl Z_{J,y_J,1}$ where $y_J$ is the longest element in
$\WW^J$ and $T_J$ is a torus acting freely on $Z_{J,y_J,1}$. Let
$a_J:Z_{J,y_J,1}@>>>G_J$ be the canonical map.

For each $J$ let $\mu_{G\T G;J}:G_J@>>>\fg^*\T\fg^*$ be the moment
map of the restriction of the $G\T G$-action on $\bG$ to $G_J$.
Let $i_J:G_J@>>>\bG$ be the inclusion. From the definitions, we
have

(a) $\L(\bG)=\sqcup_{J\sub\II}(i_J)_\bst\mu_{G\T G;J}\i(\cn^-)$
\nl and $\L(Z_{J,y_J,1})=a_J^\bst(\mu_{G\T G;J}\i(\cn^-))$. \nl
Since $\L(Z_{J,y_J,1})$ is a $T_J$-stable Lagrangian subvariety of
$T^*(Z_{J,y_J,1})$ (see 1.9), it follows that $\mu_{G\T
G;J}\i(\cn^-)$ is a Lagrangian subvariety of $T^*(G_J)$. Hence,
using 1.1(c), we see that $(i_J)_\bst\mu_{G\T G;J}\i(\cn^-)$ is a
Lagrangian subvariety of $T^*\bG$. Using this and (a) we see that

(b) {\it $\L(\bG)$ is a Lagrangian subvariety of $T^*\bG$.} \nl
From the previous proof we see that

(c) {\it The set of irreducible components of $\L(\bG)$ is in
natural bijection with \allowbreak $\sqc_{J\sub\II}\sqc_{\ss\in
S(J,\Ad(y_J))}F(\dd_{\ss_\iy})$ (notation of 1.7, 1.8).}

\subhead 1.12\endsubhead Let $X=Z_{J, y, \d}$ or $\bar{G}$. There
is a well-defined map from the irreducible components of $\L(X)$
to the nilpotent conjugacy classes of $\fg^*$ which sends the
irreducible component $C$ of $\L(X)$ to the nilpotent conjugacy
class $\co$, where $\mu_{G \times G}(C) \cap \{(f, -f) \in \co
\times \co\}$ is dense in $\mu_{G \times G}(C)$.

\enddocument